# Robust Pairwise n-Person Stochastic Duel Game

AMANG (SONG-KYOO) KIM

## Abstract

This paper introduces an extended version of a stochastic game under the antagonistic duel-type setup. The most flexible multiple person duel game is analytically solved. Moreover, the explicit formulas are solved to determine the time-dependent duel game model using the first exceed theory in multiple game stages. Unlike conventional stochastic duel games, multiple battlefields are firstly introduced and each battlefield becomes a shooting ground of pairwise players in a multi-person game. Each player selects different targets in different game stages. An analogue of this new theory was designed to find the best shooting time within multiple battlefields. This model is fully mathematically explained and is the basis with which to apply a stochastic duel-type game in various practical applications.

## 1. Introduction

Game theory has been adapted various applications in the business and management topics including economics, marketing and strategic management [1]. The domain has been expended even for the Blockchain technologies [2-6]. A conventional duel game is an arranged engagement in a combat situation between two players, with matched weapons in accordance with the agreed rules under different conditions [7]. In a versatile stochastic duel games, a joint functional of a standard stopping games has been constructed to analyze the decision-making parameters in the time domain [8] and another extended version of a multi-person stochastic duel game for multiple players is designed as another successor of the recent research by the author [9].

In this paper, the robust method has been designed for selecting best strategies in a general antagonistic n-person stochastic duel game. This general duel game is a multiple person game and each player in the game has only one bullet with two strategic choices although multiple players are joined. On his turn (or iteration), each player could choose either "wait" for one step closer to the target or "shoot" to kill one target player. A player at random times with random impacts can take the best shooting after passing the fixed threshold but players keep in mind that each player has only one chance to shoot someone.

## 2. *n*-Person Stochastic Duel Game

The antagonistic stochastic duel game of *n* players are introduced and all players know the full information regarding the success probabilities based on the time domain. A player has three strategies which are "shoot a next shooter $i_{m+1}$", "shoot a target $k_0$" and "hold". After completing each iteration, each player chooses the one within these three strategies.

**2.1 Preliminaries**

The *n*-person antagonistic duel game contains the discrete random time series (i.e., iteration period). Let $A_i(t)$ be a payoff (related) function of players $i = 1, \ldots, n$ at time $t_i$ and the payoff functions of all players are assigned as follows:

$$\{A_i(t) : 0 \leq A_i(t) \leq A_i(t+\Delta),\ t \in [0, t_i^{\max}],\ t_i^{\max} \in \mathbb{R}_+,\ \Delta > 0\}, \tag{2.1}$$

and all these functions are monotone nondecreasing. The accumulative success probability functions of players could be assigned as follows:

$$P_i(t) := \frac{A_i(t)}{A_i(t_i^{\max})},\ P_j(t) := \frac{A_j(t)}{A_j(t_j^{\max})},\ j \in -i,\ i = 1, \ldots, n. \tag{2.2}$$

The best strategy for each player in a *n* person duel game is the moment to hit others at once. As it is shown in Fig. 1, the accumulative success probabilities of all players are arbitrary incremental continuous functions which reaches 1 when the time $t_i$ meets the allowed maximum $(t_i^{\max},\ t_{-i}^{\max} < \infty)$.

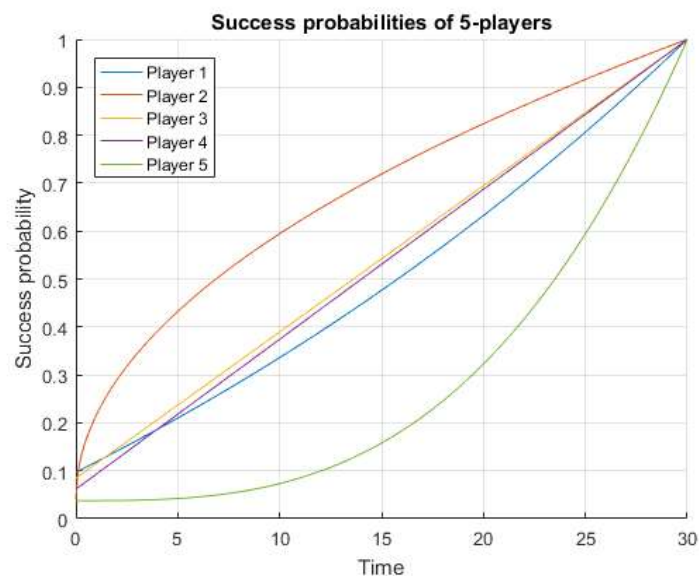

**Figure 1.** Success probabilities for a 5-person duel game

In duel games, a certain point $t^1_{\{i_1,j_1\}}$ maximizes the chance for succeeding the shoot within the first pairwise players and this best point is the first moment of the success in the continuous time domain:

$$t^1_{\{i_1,j_1\}} = \underset{t_{\{i,j\}}}{\text{argmin}} \left\{ t_{\{i,j\}} \geq 0 : \{P_{i_1}(t) - (1 - P_{j_1}(t)) \geq 0\}, j_1 \in -i_1 \right\}, \quad (2.3)$$

where $n$ is the number of players which is 5 for our case. It is noted that each player can make the decision at the certain points of the time and that is the reason why it becomes a discrete time series although the success probabilities are continuous functions.

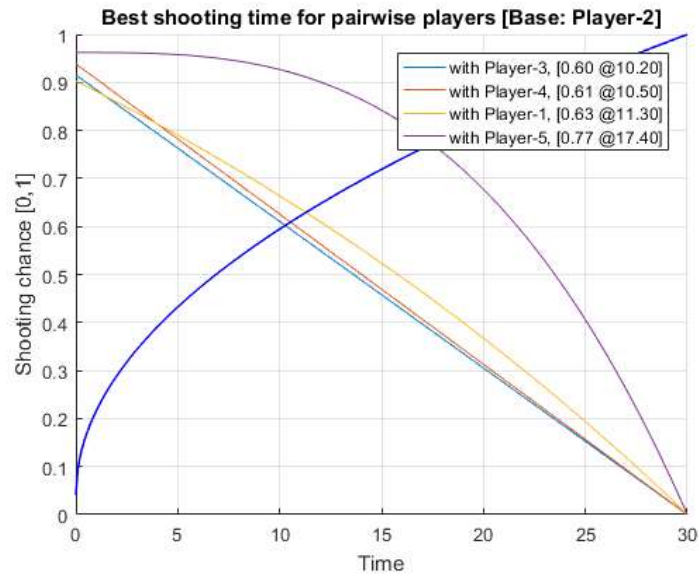

**Figure 2.** Best pairwise shooting of player 2 at multiple battle fields

In case of the 5 person duel game in Fig. 1, the first shooting pair with player 2 is player 3 and the best shooting time for player 2 is 10.2 with the 60.0 percent winning chance of the game and player 3 has 39.6 percent winning chance at same best shooting time (see Fig. 2). In a 5-person duel game, 10 battle fields (i.e., $m \in \{1, \ldots, 10\}$) for duel games are installed and each player is involved in 4 battles (or games) within 10 battles. The best shooting moments of player $i$ in the $m$-th battle field could defined as follows:

$$t^m_{\{i,j_m\}} = \underset{t_{\{i,j_m\}}}{\text{argmin}} \left\{ t_{\{i,j_m\}} \geq 0 : \{P_i(t) - (1 - P_{j_m}(t)) \geq 0\}, j_m \in -i \right\}, \quad (2.4)$$

Each players will have three more (i.e., overall $n-1$ times) chances to face duel games with other players. In the later battle field (i.e., $m > 1$), both players should be survived and hold their bullets before meeting each other. Let $(\Omega, \mathcal{F}(\Omega), P)$ be a probability time space and let $\mathcal{F}_{T^i} \subseteq \mathcal{F}(\Omega)$ be independent $\sigma$-subalgebras. Suppose

$$T^i := \sum_{k \geq 0} \varepsilon_{T^i_k}, \quad (0 =) T^i_0 < T^i_1 < \ldots, \quad i = 1, \ldots, n, \quad (2.5)$$

are $\mathcal{F}_{T^i}$-measurable renewal point processes with the following notation for players:

$$\tau^i := \begin{cases} \tau_k^i = T_k^i - T_{k-1}^i, & k = 1, 2, \ldots, \\ 0, & k \leq 0. \end{cases} \qquad (2.6)$$

We can evaluate the functional

$$\delta_i(\theta) = \mathbb{E}\left[e^{-\theta \tau^i}\right], \ \text{Re}(\theta) \geq 0, \qquad (2.6b)$$

The game in the paper is a stochastic process describing the evolution of conflicting between players with the completed information which indicates that the success at each battle field. To further formalize the game, the exit indices at the $m$-th battle field are introduced as follows:

$$\nu_m^i := \inf\{k : T_k^i = T_0^i(=0) + \tau_1 + \cdots + \tau_k^i \geq U_m^i\}, \ i = 1, \ldots, n. \qquad (2.7)$$

In the time domain of duel games, the threshold of the player $i$ could be converged into one value $t_{\{i,j_m\}}^m$ from (2.4). The player $i$ at the $m$-th battle field at $t_{\{i,j_m\}}^m$ will have the best chance to succeed for shooting compare to the failure chance of other players ($P_i\left(T_{\nu_m^i}^i\right)$ and $1 - P_{j_m}\left(T_{\nu_m^{j_m}}^{j_m}\right)$ respectively). Given the condition that both players are holding their bullets until the moment of $t_{\{i,j_m\}}^m$, player $i$ at the battle field has the highest success probability of shooting at time $T_{\nu_m^i}^i$, unless player $j_m$ does not reach his/her best shooting at the time $T_{\nu_m^{j_m}}^{j_m}$. Thus, the game is ended at $\min\left\{T_{\nu_m^i}^i, T_{\nu_m^{j_m}}^{j_m}\right\}$, $j_m \in -i$. However, we are targeting the confined duel game for player $i$ on trace $\sigma$-algebra $\mathcal{F}(\Omega) \cap \left\{P_i\left(T_{\nu_m^i}^i\right) + P_{j_m}\left(T_{\nu_m^{j_m}}^{j_m}\right) \geq 1\right\} \cap \left\{T_{\nu_m^i} \leq T_{\nu_m^{j_m}}\right\}$ (i. e., player $i$ in the game obtains the best chance for shooting first). The functional at the pairwise players at the moment of

$$\Phi_{\nu^i \nu^{-i}}^i = \Phi_{\nu^i \nu^{-i}}^i(\theta_0, \theta_1) \qquad (2.8)$$

$$= \mathbb{E}\left[e^{-\theta_0 S_{\mu-1} - \theta_1 S_\mu} \cdot \mathbf{1}_{\{T_{\nu^i} \leq T_{\nu^{-i}}\}} \mathbf{1}_{\{P_i(t) \geq 1 - P_{-i}(t)\}}\right],$$

$$Re(\theta_0) > 0, \ Re(\theta_1) > 0, \qquad (2.9)$$

where

$$\nu^i := \nu_m^i, \ \nu_m^{-i} := \{\nu_m^{j_m} : j_m \in -i_m\}, \qquad (2.10)$$

of the game will represent the status of player $i$ at the battle field $t_{\{i,j_m\}}^m$ and the rest of players upon the *exit time* $T_{\nu_m^i}^i$ and the *pre-exit time* $T_{\nu_m^i - 1}^i$ [1-4]. The Laplace-Carson transform is applied as follows:

(2.11)

$$\widehat{\mathcal{L}}_{pq}(\bullet)(u, v) = uv \int_{p=0}^{\infty} \int_{p=0}^{\infty} e^{-up - vq}(\bullet) d(p, q), \ \text{Re}(u) > 0, \ \text{Re}(v) > 0,$$

with the inverse

$$\widehat{\mathcal{L}}_{uv}^{-1}(\bullet)(p,q) = \mathcal{L}^{-1}\left(\bullet \frac{1}{uv}\right) \tag{2.12}$$

and

$$\widehat{\mathcal{L}}_{uv}^{-1}(\bullet)(r) = \widehat{\mathcal{L}}_{uv}^{-1}(\bullet)(p,q)\Big|_{(p,q)\to(r,r)} \tag{2.13}$$

where $\mathcal{L}^{-1}$ is the inverse of the bivariate Laplace transform [5]. The functional $\Phi^m_{\nu^i_m \nu^{-i}_m}$ for player $i$ at the $m$-th battle field in the game satisfies the following formula:

$$\Phi^m_{\nu^i_m \nu^{-i}_m} = \widehat{\mathcal{L}}_{uv}^{-1}\left(\mathbb{E}\left[\frac{(1-\gamma_0(\tau))(1-\Gamma_0(\sigma))}{\gamma_1(\tau)(1-\gamma_2(\tau))\Gamma_1(\sigma)(1-\Gamma(\sigma))}\right]\Gamma(t^i_m)\right)(t^i_m), \tag{2.14}$$

where

$$\left(t^{j_m}_m = \right)t^i_m := t^{\{i,j_m\}}_m \tag{2.15}$$

$$\gamma(x,t) = e^{-xt}, \tag{2.16}$$
$$\gamma(t) := \gamma(v,t), \tag{2.17}$$

$$\Gamma_0(t) := \gamma(u,t), \tag{2.18}$$
$$\Gamma_1(t) := \gamma(\theta_0 + u, t), \tag{2.19}$$
$$\Gamma_2(t) := \gamma(\theta_0 + \theta_1 + u, t), \tag{2.20}$$
$$\Gamma(t) := \gamma(t) \cdot \Gamma_2(t). \tag{2.21}$$

Introduce the families:

$$\nu^i_m(p) = \inf\{l : T^i_l > p\}, \tag{2.22}$$

$$\nu^{-i}_m(q) = \{k : T^{j_m}_k > q\}, \tag{2.23}$$

The functional $\Phi^i_{\nu^i \nu^{-i}}$ for the player set $\{i_m, j_m(=-i)\}$ has all decision making parameters of player $i_m$. The information includes the best shooting moments ($T^i_{\nu^i}$; *exit time*), the one step prior to each best shooting moment ($T^i_{\nu^i_m - 1}$; *pre-exit time*) and the optimal number of iterations for player $i$ at $m$-th battle field. The information from the closed functional are as follows:

$$\mathbb{E}\left[e^{-\theta T^i_{\nu^i_m}}\right] = \Phi^m_{\nu^i_m \nu^{-i}_m}(0,\theta), \tag{2.24}$$

$$\mathbb{E}\left[T^i_{\nu^i_m}\right] = \lim_{\theta \to 0}\left(-\frac{\partial}{\partial \theta}\right)\Phi^m_{\nu^i \nu^{-i}}(0,\theta), \tag{2.25}$$

$$\mathbb{E}\left[T^i_{\nu^i_m - 1}\right] = \lim_{\theta \to 0}\left(-\frac{\partial}{\partial \theta}\right)\Phi^m_{\nu^i \nu^{-i}}(\theta, 0). \tag{2.26}$$

### 2.2. Pairwise Set of Shooting Orders

Although this game is designed for multiple players, the winning chance of each player is independent with other players. Each player has his/her own best chance to win a game because of the backward induction in any duel type games [1-5]. Hence, finding the best shoot moment when the sum of success probability of a player passes the threshold is the most critical matter [2] regardless of his/her shooting performance. The best shooting moments of the players are determined from (2.4) and the sorted list of pairwise player set at the beginning is as follows:

$$M^0 = \underset{\{i_m, j_m\} \in N}{\text{argsort}} \left\{ \{i_m, j_m\}_m \middle| t^1_{\{i_1,j_1\}} \leq \cdots \leq t^m_{\{i_m,j_m\}} \cdots \leq t^{m_0}_{\{i_{m_0},j_{m_0}\}} \right\}, \quad (2.27)$$

$$t^0 = \left\{ t^1_{\{i_1,j_1\}} \leq \cdots \leq t^m_{\{i_m,j_m\}} \cdots \leq t^{m_0}_{\{i_{m_0},j_{m_0}\}} \right\}, \quad (2.28)$$

where $m_0 = \binom{n}{2}$ and $i_m$ is the $m$-th shooting player. The shooting orders of one specific player $i$ could be similarly defined as follows:

$$M^0_i = \underset{\{i, r\}}{\text{argsort}} \left\{ ^\exists \{i, r\}_{m_r} \middle| t^{m_r}_{\{i,r\}}, \ r = -i, \ m_{r \in \{1,2,\ldots,n-1\}} \right\}, \quad (2.29)$$

$$t^0_i = \underset{t_{\{i,r\}}}{\text{argsort}} \left\{ ^\exists t^{m_r}_{\{i,r\}} \middle| r = -i, \ m_{r \in \{1,2,\ldots,n-1\}} \right\}, \quad (2.30)$$

and the numbers of the set $M^0_i$ and $t^0_i$ are $n-1$ (i.e., $n(M^0_i) = n(t^0_i) = n-1$). The shooting order of the demonstrated set of players (see Fig. 1) is shown in Fig. 4.

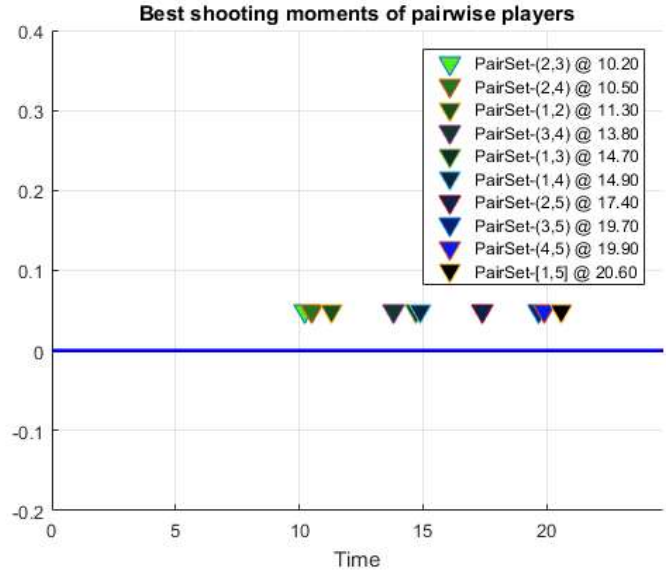

**Figure 4.** The shooting order of a 5-person duel game

As it is illustrated on Fig. 4, the first shooter are either player 2 or player 3 whoever passes the first best shooting moment ($t^1_{\{2,3\}} = 10.2$) first. The next pairwise set of best shooters is the player set $\{2, 4\}$ if the first optimal shooting moment passes without any shooting from the first pairwise players. If either player 2 or player 3 shoots at the first shooting moment, the whole player status of the game shall be changed. From (2.4) and (2.27), the sorted shooting list of all player set is follows:

$$\boldsymbol{M^0} = \{\{2,3\}_1, \cdots, \{i_m, j_m\}_m, \cdots, \{1,5\}_{10}\}, \ i_\bullet \neq j_\bullet \qquad (2.31)$$

and the moments of the best shooting for all players are as follows:

$$\boldsymbol{t^0} = \left\{ t^1_{\{2,3\}} \leq t^2_{\{2,4\}} \leq \cdots \leq t^{10}_{\{1,5\}} \right\} \qquad (2.32)$$

$$= \{10.2 \leq 10.5 \leq 11.30 \leq \ldots \leq 20.50\},$$

and these best shooting moments are equivalent with the moments when the battle fields are opened (i.e., the $m$-th battle is opened at the moment $t^m_{\{\bullet\}}$). In the case of player $i$ (the player 3 in our demonstration), the shooting list of this player could be determined as follows from (2.29)-(2.30):

$$\boldsymbol{M^0_2} = \{\{2,3\}_1, \{2,4\}_2, \{2,1\}_3, \{2,5\}_7\}. \qquad (2.33)$$

and

$$\boldsymbol{t^0_2} = \left\{ t^1_{\{2,3\}} \leq t^2_{\{2,4\}} \leq t^3_{\{2,1\}} \leq t^7_{\{2,5\}} \right\} \qquad (2.34)$$

$$= \{10.2 \leq 10.5 \leq 11.3 \leq 17.4\}.$$

The details about the strategies of each player are covered on Section 3.

**2.3. Status change of pairwise $n$-person duel game after shooting**

Let us consider the $n$-person antagonistic duel game with the pure strategy. Although this duel game is targeted for multiple players, the pairwise two-person game between players $i_m$ and $j_m$ is considered from the shooting order. The set of $n$ players which is sorted by the shooting order at the beginning could be defined as follows:

$$\boldsymbol{M^0} = \{\{i_1, j_1\}_1, \ldots, \{i_m, j_m\}_m, \ldots \{i_{m_0}, j_{m_0}\}_{m_0}\}. \qquad (2.35)$$

where

$$m_0 = 1, \ldots, \binom{n}{2}. \qquad (2.36)$$

Unlike an OneToN game [4], each player has only one bullet which makes a duel game more complicated. Hence, the status of players are changed after shooting. For instance, if the player $i_m$ shoots the target $j_m$, the player status is changed and the sorted player set is updated as follows:

$$(\boldsymbol{N} \setminus \{j_m\} =) \boldsymbol{N_1} = \{1, .., j_m - 1, j_m + 1, \ldots, i_n\}, \qquad (2.37)$$

$$t^1_{\{i_1,j_1\}} \leq \cdots \leq t^w_{\{i_w,j_w\}} \cdots \leq t^{w_0}_{\{i_{w_0},j_{w_0}\}}, \ w_0 = m_0 - 2. \qquad (2.38)$$

It is noted that the success probability shall be changed to zero when a bullet is exhausted. We could summarized the rules of this duel type game as follows:

1. Each player has a gun with single bullet;
2. The player $i_1$ who is the first ordered shooter at the first battlefield;

3. Each player could have two strategic choices;
4. Once shooting is occurred, the success probability of a shooter shall be zero;
5. The game ends when all players use their bullet;

According to the above rules, the hitting probabilities of both players are recursive for every status change. Hence, the success probabilities of all players are reset and starting at the beginning (i.e., $t \in [0, t_i^{\max}]$). Previously mentioned, the success status of the game shall be totally changed once a player $i_m$ shoots the target (i.e., $j_m$). Let us consider the set of success probability matrix $\mathbb{P}(t, \boldsymbol{N})$ as follows:

$$\mathbb{P}(t, \boldsymbol{N}) = \{P_{i_1}(t), P_{i_2}(t), P_{i_m}(t), \ldots, P_{i_n}(t)\}. \tag{2.39}$$

When player $i_m$ shoots player $j_m$, $\mathbb{P}(t, \boldsymbol{N})$ shall be changed as follows:

$$\mathbb{P}(t, \boldsymbol{N_1}) = \begin{cases} \{P_{i_1}(t), \ldots, \boldsymbol{0}_{i_m}(t), \ldots, P_{i_n}(t)\}, & \textit{fail}, \\ \mathbb{P}(t, \boldsymbol{N} \setminus \{j_m\}), & \textit{succeed}. \end{cases} \tag{2.40}$$

and the pairwise player set after shooting has been changed as follows:

$$\boldsymbol{M^1} = \begin{cases} \boldsymbol{M^0}, & \textit{fail}, \\ \boldsymbol{M^0} \setminus \{\{j, k_m\}_m, k_m = 1, \ldots, n-2, \ k_m \neq j\}, & \textit{succeed}. \end{cases} \tag{2.41}$$

## 3. Optimal Strategies for $n$-Person Duel Games

Unlike to analyze a conventional multi-person game, a duel type multi-person game could analyze the best strategies for each player in the early stage and it could be described as a conventional two-person game is intuitively easy to visualize its game system. Because the shooting orders of all players are determined by backward induction (see Section 2.2), the pairwise player set for each battle field and the best shooting times with its competitors in each pairwise player set are determined from (2.29)-(2.30).

### 3.1. Finding the best battle field for shooting

Each player has $n-1$ battle field at the best shooting times and faces different opponent players. Recalling from (2.30), the set of shooting order for player $i$ is constructed as follows:

$$\boldsymbol{t_i^0} = \left\{ t_{\{i,j_{m_1}\}}^{m_1} \leq \cdots \leq t_{\{i,j_{m_r}\}}^{m_r} \leq \cdots \leq t_{\{i,j_{m_{n-1}}\}}^{m_{n-1}} \right\}. \tag{3.6}$$

In a two-person duel type game [8,9], the probability of successful shooting is equivalent with the complement of the failing probability of an opponent player and this same concept is applied in a pairwise player duel game at the $m$-th battle field. Hence, the actual successful shooting probability of a player $i$ (with an opponent player $j$) at the $m$-th battle field could be found as follows:

$$\mathcal{P}_{ij}^m = 1 - P_j\left(t_{\{i,j\}}^m\right) \left\{ \prod_{h=1}^{m-1} P_j\left(t_{\{i_h,j_h\}}^h\right)^{\mathbf{1}_{\{j=j_h\}}} \right\}, \tag{3.7}$$

and the best battle field indicator as follows:

$$\varrho_{ij}^m = \frac{\mathcal{P}_{ij}^m}{\mathcal{P}_{ji}^m}, \; j = -i, \; m = 1, 2, \ldots, \binom{n}{2}, \qquad (3.8)$$

where $n$ is the number of players. From (3.6)-(3.8), the best battle field and the best moment of shooting for player $i$ is as follow:

$$m_i^* = \underset{m}{\operatorname{argmin}}\left\{m \,\middle|\, \varrho_{ij}^m, \; j = -i\right\}, \; t_i^* = t_{\{i, j_{m_i^*}\}}^{m_i^*}, \qquad (3.9)$$

and the best target for shooting is the player $j_{m_i^*}$ at the $m_i^*$-th battle field. The main implication is that a player selects the target who have the highest chance for successful shooting him within his opponent players. It is noted that the situation for multiple bullets ($b < n - 1$) is easily extended. From (3.8), the set of battle field for shooting range could be determined as follows:

$$\boldsymbol{m}_i^* = \underset{m}{\operatorname{argsort}}\left\{m \,\middle|\, \varrho_{ij}^{m_1} \leq \ldots \leq \varrho_{ij}^{m_b}, \; j = -i, r = 1, \ldots, b\right\}, \qquad (3.10)$$

where $b$ is the total number of available bullets for player $i$. It is noted that player $i$ can shoot at ever battle field if he has more bullets than joined players. Although this analysis provides the optimal decision parameter, it is only applicable for best moment, not for the final result. The actual result of the game still requires additional analysis based on the two-person versatile duel game [8] which has been previously analyzed on the next section.

**3.2. Best shooting moment in the $m$-th battle field**

Once a battle field is fixed from Section 3.1, the functional $\varPhi_{\nu^i\nu^{-i}}^i$ gives the full analytical information to build up the winning strategies for pairwise players in the battle field. Recalling from (2.24)-(2.26), the functional $\varPhi_{\nu^i\nu^{-i}}^i$ for the player set $\{i, j(=-i)\}$ in the $m$-th battle field (i. e., the subscript character $m$ in each notation is dropped) has all decision making parameters of player $i_m$. The information includes the best moments of shooting ($T_{\nu^i}^i$, $T_{\nu^{-i}}^j$; *exit time*), the one step before the best moment*s* ($T_{\nu^i-1}^i$, $T_{\nu^{-i}-1}^j$; *pre-exit time*) and the optimal number of iterations for both players. The information for both players from the closed functional are as follows:

(for player $i$ at the $m$-th battle field )

$$\mathbb{E}\left[T_{\nu^i}^i\right] = \lim_{\theta \to 0}\left(-\tfrac{\partial}{\partial \theta}\right)\varPhi_{\nu^i\nu^{-i}}(0, \theta, 0, 0), \qquad (3.1)$$

$$\mathbb{E}\left[T_{\nu^i-1}^i\right] = \lim_{\theta \to 0}\left(-\tfrac{\partial}{\partial \theta}\right)\varPhi_{\nu^i\nu^{-i}}(\theta, 0, 0, 0), \qquad (3.2)$$

(for player $j(=-i_m)$ at the $m$-battle field )

$$\mathbb{E}\left[T_{\nu^j}^j\right] = \lim_{\vartheta \to 0}\left(-\tfrac{\partial}{\partial \theta}\right)\varPhi_{\varPhi_{\nu^i\nu^{-i}}}(0, 0, 0, \vartheta), \qquad (3.3)$$

$$\mathbb{E}\left[T_{\nu^j-1}^j\right] = \lim_{\vartheta \to 0}\left(-\tfrac{\partial}{\partial \theta}\right)\varPhi_{\varPhi_{\nu^i\nu^{-i}}}(0, 0, \vartheta, 0), \qquad (3.4)$$

and

$$\mathbb{E}\left[\nu^i\right] \simeq \left\lfloor \mathbb{E}\left[\frac{T_{\nu^i}^i}{\tau^i}\right] \right\rfloor, \mathbb{E}\left[\nu^j\right] \simeq \left\lfloor \mathbb{E}\left[\frac{T_{\nu^j}^j}{\tau^j}\right] \right\rfloor. \tag{3.5}$$

According to a conventional two-person duel game strategies, player $i$ should take a shoot at his $(\nu^i)$-th turn (which means he should wait until the $(\nu^i - 1)$-th iterations). The average duration for player $i$ to get the best shooting chance at the $m$-th battle field becomes $\mathbb{E}\left[T_{\nu^i}^i\right]$. Similarly, player $j$ should shoot at his $(\nu^j)$-th iteration and the average duration is $\mathbb{E}\left[T_{\nu^j}^j\right]$. Each player will continue their iterations until the accumulated iteration passes the threshold moment $t_{\{i,j\}}^m$ [7]. Additionally, player $i$ should compare the chances not only at $T_{\nu^i}^i$ but also at $T_{\nu^i-1}^i$ to build up a proper strategy. In the adverse condition of player $i$ at the $m$-th battle field (i.e., $\left(t_{\{i,j\}}^m \leq \right)\mathbb{E}\left[T_{\nu^i}^i\right] < \mathbb{E}\left[T_{\nu^j}^j\right]$), player $i$ shall take the shoot at $\mathbb{E}\left[T_{\nu^i}^i\right]$ after player $j$ fails his shoot at $\mathbb{E}\left[T_{\nu^j}^j\right]$ if player $i$ has higher chance at $\mathbb{E}\left[T_{\nu^i}^i\right]$ compare to $\mathbb{E}\left[T_{\nu^i-1}^i\right]$ (i.e., $\mathbb{E}\left[P_i\left(T_{\nu^i-1}^i\right)\right] < \mathbb{E}\left[P_i\left(T_{\nu^i}^i\right)\right]$). Otherwise, player $i$ should take a shoot at the time $\mathbb{E}\left[P_i\left(T_{\nu^i-1}^i\right)\right]$ instead of the time $\mathbb{E}\left[T_{\nu^i}^i\right]$ (i.e., $\mathbb{E}\left[P_i\left(T_{\nu^i-1}^i\right)\right] \geq \mathbb{E}\left[P_i\left(T_{\nu^i}^i\right)\right]$) [7]. Because of the backward induction, it does not matter players are better or worse shooting but the matter is the cumulative success probabilities [7].

## ACKNOWLEDGEMENT
The corresponding Matlab codes are publicly available on the GitHub (https://github.com/amangkim/rpnsdg) for users to try the demonstrations of the robust pairwise n-person stochastic duel game.

## REFERENCES
[1] Moschini, G. Nash equilibrium in strictly competitive games: Live play in soccer. Econ. Lett. 2004, 85, 365-371.
[2] Johnson, B., Laszka, A. and et al., Game-theoretic analysis of ddos attacks against bitcoin mining pools, International Conference on Financial Cryptography and Data Security. Springer, 2014, pp. 72-86.
[3] Laszka, A., Johnson, B. and Grossklags, J., When bitcoin mining pools run dry, International Conference on Financial Cryptography and Data Security. Springer, 2015, pp. 63-77.
[4] Kim, S. -K., Yeun, C. Y. and et al., Blockchain Adoption For Automotive Security by Using Systematic Innovation, IEEE Proceedings of ITEC Asia-Pacific 2019, pp. 1-4.
[5] Kim, S.-K., Blockchain Governance Game, Comp. Indust. Eng., 2019, 136, pp. 373-380.
[6] Kim, S.-K., Strategic Alliance For Blockchain Governance Game, Probability in the Engineering and Informational Sciences, 2020, Accepted.
[7] Polak, B. Discussion of Duel. Open Yale Courses. 2008. Available online: http://oyc.yale.edu/economics/econ-159/lecture-16 (accessed on 1 May 2019).
[8] Kim, S.-K. A Versatile Stochastic Duel Game, Mathematics, 2020, 8, 678.
[9] Kim, S.-K. Antagonistic One-To-N Stochastic Duel Game, Mathematics, 2020, 8, 1114.